\documentstyle[12pt]{article}
\title{On Braid Words and Irreflexivity \\
\author{{\sc David M.~Larue} \\
       Department of Mathematics \\
       University of Colorado at Boulder \\
       Boulder, Colorado  80309-0426  USA}
\date{} % 4-7-92
}
\newcommand{\bproof}{\hspace*{-\parindent}$\vdash$}
\newcommand{\eproof}{\hspace*{\fill}$\dashv$}
\newtheorem{lemma}{Lemma}

\newtheorem{definition}[lemma]{Definition}
\newtheorem{theorem}[lemma]{Theorem}
\newtheorem{case}{Case}[lemma]
\newtheorem{subcase}{SubCase}[case]
\newcommand{\shift}{s}
\newcommand{\bid}{\varepsilon}
\newcommand{\bop}{\cdot}
\newcommand{\bgp}{B_\infty}
\newcommand{\bel}{\sigma}
\newcommand{\fid}{\varepsilon}
\newcommand{\fop}{\cdot}
\newcommand{\fgp}{{\cal F}_G}
\newcommand{\fel}{g}
\newcommand{\fgn}{G}
\newcommand{\xid}{\varepsilon}
\newcommand{\xop}{\cdot}
\newcommand{\xgp}{{\cal F}_X}
\newcommand{\xel}{x}
\newcommand{\xgn}{X}
\begin{document}
\maketitle

A left-distributive algebra is a set ${\cal B}$ equipped with a 
binary operation (here written as concatenation)
such that $a(bc)=(ab)(ac)$ for all $a,b,c\in{\cal B}$.  The free
left-distributive algebra on $n$ letters is denoted ${\cal A}_n$,
and we write ${\cal A}$ for ${\cal A}_1$

If $P,Q\in{\cal B}$ write $P<_L Q$ iff one can write $P$ as
a strict prefix of $Q$, i.e., $Q=((PQ_1)\ldots)Q_k$ for some
$Q_1,\ldots,Q_k,k\ge1$.  Then a proof that $<_L$ is irreflexive
(that is, that $P\ne ((PQ_1)\ldots)Q_k$ for all $P,Q_1,\ldots,
Q_k\in{\cal A}$) on ${\cal A}$ was found by
R.~Laver (\cite{lav:ld}), under large cardinal assumptions, as part of
a theorem that ${\cal A}$ is isomorphic to a certain algebra of
elementary embeddings from set theory.

It was also proved in \cite{lav:ld} that $<_L$ linearly orders 
${\cal A}$, the part that for all $P,Q\in{\cal A}$ at least one of
$P<_L Q,P=Q,Q<_L P$ holds being proved independently and by a
different method by P.~Dehornoy (\cite{deh:fdg,deh:gl}).  The linear
ordering of ${\cal A}$ gives left cancellation, the solvability of the
word problem, and other consequences.  Left open was whether
irreflexivity, and hence the linear ordering, can be proved in ZFC.

Recently, Dehornoy (\cite{deh:bg}) has found such a proof, 
involving an extension of the infinite braid group but without
invoking axioms extending ZFC.  
The purpose of this note is to prove the result without the
additional machinery of this extended group, and at shorter length. 

\section{The $\sigma_1$-proposition implies irreflexivity}

We recall from \cite{deh:bg} the connection between the braid group
and ${\cal A}$.

Note that since ${\cal A}$ is free, 
if for some left-distributive algebra ${\cal B}$,
$<_L$ is irreflexive on ${\cal B}$, then $<_L$ is irreflexive on
${\cal A}$.  Dehornoy found such a ${\cal B}$, where ${\cal B}$ is a
subset of the infinite dimensional braid group endowed with a bracket
operation.

\begin{definition}
The infinite braid group $\langle \bgp,\bid,\bop\rangle$ is given 
by generators $\{\bel_i:i\in\omega^+\}$ and relations 
$\bel_i\bop\sigma_j=\bel_j\bop\bel_i$ when $|i-j|>1$ and 
$\bel_i\bop\bel_j\bop\bel_i= \bel_j\bop\bel_i\bop\bel_j$ when $|i- 
j|=1$.
\end{definition}

\begin{definition}
$\shift$ is the endomorphism of $\bgp$ which extends the shift map
$\shift(\bel_i)=\bel_{i+1}$ (possible because $\shift$ preserves
the defining relations of $\bgp$).
\end{definition}

\begin{definition}
(The Dehornoy bracket)  For $p,q\in\bgp$ define a binary operation
$[\,]$: \[p[q]=p\bop\shift(q)\bop\bel_1\bop\shift(p)^{-1}\].
\end{definition}

Motivation:  Suppose that ${\cal B}$ is a left-distributive algebra
with left cancellability (as ${\cal A}$ turns out to be from the
linearity of $<_L$).  Then the braid group has a partial action on
${\cal B}^\omega$ defined by, for $\vec{b}\in{\cal B}^\omega$,
\[
  \left((\vec{b})\sigma_i\right)_j =
\left\{
  \begin{array}{ll}
    \vec{b}_i\vec{b}_{i+1} & {\rm\ if\ }j=i \\
    \vec{b}_i         & {\rm\ if\ }j=i+1 \\
    \vec{b}_j         & {\rm\ else.}
  \end{array}
\right.
\]
with left cancellability making the partial actions of the inverses of braid
generators well defined.
Thus $(P,Q,R,S,\ldots)\sigma_2=(P,QR,Q,S,\ldots)$
for $P,Q,R,S,\ldots\in{\cal B}$.

Then one is led to the Dehornoy bracket in the braid group 
by noting that, for $x,P,Q\in{\cal B}$ and $p,q\in\bgp$, 
if $(x,x,x,\ldots)p=(P,x,x,\ldots)$ and 
$(x,x,x,\ldots)q=(Q,x,x,\ldots)$ then
\begin{eqnarray*}
  (x,x,x,\ldots)p[q] &=& 
(x,x,x,\ldots)p\bop\shift(q)\bop\bel_1\bop\shift(p)^{-1} \\
 &=& (P,x,x,\ldots)\shift(q)\bop\bel_1\bop\shift(p)^{-1} \\
 &=& (P,Q,x,\ldots)\bel_1\bop\shift(p)^{-1} \\
 &=& (PQ,P,x,\ldots)\shift(p)^{-1} \\
 &=& (PQ,x,x,\ldots)
\end{eqnarray*}
which suggests, in order that $(x,x,x,\ldots)p[q]=(PQ,x,x,\ldots)$,
the above definition of $p[q]$.

Dehornoy proved as a corollary of his irreflexivity result
the following theorem:

\begin{theorem}\label{sigma}
(Dehornoy---the $\sigma_n$-proposition) If $p\in\bgp$ is written as a 
product of generators and their inverses, including at least one 
$\bel_n$ and no $\bel_n^{-1}$, then $p\ne\bid$.
\end{theorem}

He also showed that a direct proof of theorem \ref{sigma} would give
a relatively short proof of irreflexivity, namely:

\begin{theorem}\label{irr2}
(Dehornoy)
\begin{enumerate}
\item The Dehornoy bracket $[\,]$ on the braid group $\bgp$ is 
left distributive;
\item If the $\sigma_1$-proposition holds then irreflexivity holds in
$\bgp$ and so also in ${\cal A}$;
furthermore, the closure under the Dehornoy bracket of any 
element of the braid group is isomorphic to ${\cal A}$.
\end{enumerate}
\end{theorem}

\bproof
Proof of \ref{irr2}.1:
Compute
\[p[q[r]] =p\bop\shift(q)\bop\shift^2(r)\bop\bel_2\bop\bel_1
     \bop\shift^2(q)^{-1}\bop\shift(p)^{-1}
=p[q][p[r]] \]
to get left-distributivity.

Proof of \ref{irr2}.2:
Suppose $p,q_1,\ldots,q_k\in\bgp$ satisfy
$p=p[q_1]\ldots[q_k]$.  The right hand expands to
$p\bop(r_1\bop\bel_1\bop s_1)\bop(r_2\bop\bel_1\bop s_2)\bop\ldots\bop
(r_k\bop\bel_1\bop s_k)$ with $r_i,s_i\in\shift(\bgp)$.
Then multiplying by $p^{-1}$ one obtains a contradiction to the
$\sigma_1$-proposition.  Therefore irreflexivity holds in $\bgp$, and
hence in ${\cal A}$.

Let $x$ be the generator of the free left-distributive algebra ${\cal
A}$.  For $r\in\bgp$ define the homomorphism $\chi^r:{\cal
A}\rightarrow\langle\bgp,[\,]\rangle$ inductively by $\chi^r_x=r$ and
$\chi^r_{PQ}=\chi^r_P[\chi^r_Q]$.  This is well defined as $[\,]$ is
left-distributive and ${\cal A}$ is free.  To show that $\chi^r$ is
injective, if $P\ne Q$ in ${\cal A}$, by trichotomy suppose $P<_L
Q$.  Conclude $\chi^r_P<_L \chi^r_Q$, and so by irreflexivity in
$\bgp$ that $\chi^r_P\ne\chi^r_Q$.  Hence $\langle r,[\,]\rangle$ is
isomorphic to ${\cal A}$.
\eproof

\section{Proof of the $\sigma_1$-proposition}

We will prove the $\sigma_1$-proposition first, using the action of
the braid group $\bgp$ as automorphisms of the free group on
countably many generators, getting as quickly as possible to the
minimum needed for irreflexivity.  Subsequently, we will prove the
full $\sigma_n$-proposition, for which we use a somewhat non-standard
action, and have to deal with a greater number of cases.

\begin{definition}
$\langle \fgp,\fid,\fop\rangle$ is the free group on generators 
$\fgn=\{\fel_i:i\in\omega\}$.
\end{definition}

\begin{definition}\label{fgpact}
\begin{eqnarray*}
  (\fel_i)\bel_i&=&\fel_i\fop \fel_{i+1}\fop \fel_i^{-1} \\
  (\fel_{i+1})\bel_i&=&\fel_i \\
  (\fel_j)\bel_i&=&\fel_j{\rm\ if\ }j\ne i,i+1.
\end{eqnarray*}
\end{definition}

The following is well known, see for example \cite[cor.~1.8.3,
pg.~25]{bir:bg} (although our group has an extra generator $\fel_0$
for convenience later, and we do not require faithfulness).

\begin{lemma}
The action of the $\bel_i$s on the $\fel_j$s extends to a faithful
action of $\bgp$ on $\fgp$.
\end{lemma}

We record for convenience the action of $\bel_i^{-1}$:
\begin{eqnarray*}
  (\fel_i)\bel_i^{-1}&=&\fel_{i+1} \\
  (\fel_{i+1})\bel_i^{-1}&=&\fel_{i+1}^{-1}\fop\fel_i\fop\fel_{i+1} \\
  (\fel_j)\bel_i^{-1}&=&\fel_j{\rm\ if\ }j\ne i,i+1.
\end{eqnarray*}

In passing we observe without proof that if $\xi$ is the
anti\-auto\-mor\-phism of $\bgp$ obtained by reversing products of
generators and their inverses (so
$\xi(\bel_1\bop\bel_2^{-1})=\bel_2^{-1}\bop\bel_1$), then
\[(\fel_1,\fel_2,\ldots)\xi(p)=((\fel_1)p,(\fel_2)p,\ldots),\] where
the left is the action previously defined of a braid on an element of
$\fgp^\omega$ where $\fgp$ is equipped with the left distributive
operation of conjugation ${}^{f_1}f_2=f_1\fop f_2\fop f_1^{-1}$, and
the right are the actions of braids on elements of $\fgp$ just
defined in \ref{fgpact}.

The main property about this action that we need is the following
observation:

\begin{lemma}\label{tech1}
If a reduced word $f$ in $\fgp$ begins with $\fel_1$, and
$\sigma\in\bgp$ is a generator or the inverse of a generator, but not
$\bel_1^{-1}$, then the reduced form of $(f)\sigma$ also begins with
$\fel_1$.
\end{lemma}

\bproof
Assume that the action of $\sigma$ is applied to each generator or
its inverse in the reduced form of $f$, and then a fixed reduction
is applied to the term that results to produce the reduced form of
$(f)\sigma$.

There are two cases where the reduced form of $(f)\sigma$ may fail to
begin with $\fel_1$:

\begin{case}
$\sigma=\bel_i^{\pm1}$ with $i>1$.
\end{case}

In this case, all $\fel_1^{\pm1}$s are unchanged by the action of
$\sigma$ on $f$, and none are produced.  Therefore the only way that
the leading $\fel_1$ could be cancelled in the reduction of
$(f)\sigma$ is by a $\fel_1^{-1}$ already present in the reduced
form of $f$.

Let $f=\fel_1\fop f_1\fop\fel_1^{-1}\fop f_2$, in reduced form,
display that instance of $\fel_1^{-1}$.  Then 
$(f)\sigma = \fel_1\fop(f_1)\sigma\fop\fel_1^{-1}\fop(f_2)\sigma$,
so $(f_1)\sigma=\fid$.  But then $f_1=\fid$, so that $f$ was not in
reduced form, contradiction.

\begin{case}
$\sigma=\bel_1$.
\end{case}

In this case, $\fel_1^{-1}$s are produced by
$(\fel_1^{\pm1})\bel_1=\fel_1\fop\fel_2^{\pm1}\fop\fel_1^{-1}$ and
by $(\fel_2^{-1})\bel_1=\fel_1^{-1}$.

There are two possibilities for the leading $\fel_1$ produced by the
action of $\bel_1$ on $f$ to be cancelled in the reduction:

\begin{subcase}
The leading $\fel_1$ in the unreduced form of $(f)\bel_1$ is
cancelled by a $\fel_1^{-1}$ produced by the action of $\bel_1$ on
a $\fel_1^{\pm1}$.
\end{subcase}

First we note that the $\fel_1^{\pm1}$ occurrence of $f$ which produces the
$\fel_1^{-1}$ which cancels the leading $\fel_1$ of $(f)\bel_1$ 
is distinct from the leading $\fel_1$ of $f$.

As above, let $f=\fel_1\fop f_1\fop\fel_1^{\pm1}\fop f_2$, in reduced
form, display that instance of $\fel_1^{\pm1}$.
Then
$(f)\bel_1 = \fel_1\fop\fel_2\fop\fel_1^{-1}\fop(f_1)\bel_1
 \fop \fel_1\fop\fel_2^{\pm1}\fop\fel_1^{-1}\fop(f_2)\bel_1$,
so 
$\fel_2\fop\fel_1^{-1}\fop(f_1)\bel_1
 \fop \fel_1\fop\fel_2^{\pm1}=\fid$.  But then
$f_1=\fel_1^{-1}\fop\fel_1^{\mp1}$, so that
$f$ was not in reduced form, contradiction.

\begin{subcase}
The leading $\fel_1$ in the unreduced form of $(f)\bel_1$ is
cancelled by a $\fel_1^{-1}$ produced by the action of $\bel_1$ on
a $\fel_2^{-1}$.
\end{subcase}

Again, let $f=\fel_1\fop f_1\fop\fel_2^{-1}\fop f_2$, in reduced
form, display that instance of $\fel_2^{-1}$.
Then
$(f)\bel_1 = 
\fel_1\fop\fel_2\fop\fel_1^{-1}
\fop(f_1)\bel_1\fop\fel_1^{-1}\fop(f_2)\bel_1$,
so 
$\fel_2\fop\fel_1^{-1}\fop(f_1)\bel_1=\fid$.  But then
$f_1=\fel_1^{-1}\fop\fel_2$, so that
$f$ was not in reduced form, contradiction.
\eproof

\bproof
To prove the $\sigma_1$-proposition, 
assume that ${p}\in\bgp$ is formed from a product of generators of 
$\bgp$ and their inverses, with at least one $\bel_1$ and no 
$\bel_1^{-1}$.
Write ${p}={p}_1\bop\bel_1\bop{p}_2$ 
where $\bel_1^{\pm1}$ doesn't occur in $p_1$ and $\bel_1^{-1}$
doesn't occur in $p_2$.

Then
\begin{eqnarray*}
 (\fel_1^{-1}){p} &=& (\fel_1^{-1}){p}_1\bop\bel_1\bop{p}_2 \\
  &=& (\fel_1^{-1})\bel_1\bop{p}_2 \\
  &=& (\fel_1\fop\fel_2^{-1}\fop\fel_1^{-1}){p}_2
\end{eqnarray*}

Since $\fel_1\fop\fel_2^{-1}\fop\fel_1^{-1}$ is in reduced form and
begins with $\fel_1$, and since $\bel_1^{-1}$ does not appear in
${p}_2$, then by repeated application of lemma \ref{tech1}
it follows that the reduced form of
$(\fel_1\fop\fel_2^{-1}\fop\fel_1^{-1}){p}_2$, and hence of 
$(\fel_1^{-1}){p}$, must begin with $\fel_1$.  Therefore
$(\fel_1^{-1}){p}\ne\fel_1^{-1}$, and so ${p}\ne\bid$.
\eproof

\section{Proof of the $\sigma_n$-proposition}

We now prove the full $\sigma_n$ proposition.  It can be done by the
braid action defined above, using as an invariant that the reduced
form begins with $\prod_{i=0}^n\fel_i\fop\fel_m^{\pm1}$ for some
$m>n$, and showing that words retains that property when acted upon
by a braid generator or its inverse, except for $\bel_n^{-1}$,
obtaining the desired result by starting with $\prod_{i=0}^n\fel_i$. 

We make use of a slightly different group action which decreases the
number of cases to check.  The map $\phi:\xgp\rightarrow\fgp$ given
by $(\xel_i)\phi=\prod_{n=0}^i\fel_n$, with domain the free group
$\langle \xgp,\xid,\xop\rangle$ on generators 
$\xgn=\{\xel_i:i\in\omega\}$, extends to an isomorphism with inverse
$(\fel_0)\phi^{-1}=\xel_0$,
$(\fel_{i+1})\phi^{-1}=\xel_i^{-1}\xop\xel_{i +1}$.

Then the action of $\bgp$ on $\fgp$ induces an action on $\xgp$
by, for $f\in\xgp,p\in\bgp$, defining $(f)p=(f)\phi p\phi^{-1}$.

We record the actions of generators and their inverses of $\bgp$ 
on generators of $\xgp$:
\begin{eqnarray*}
  (\xel_i)\bel_i^{\pm1}&=&\xel_{i\pm1}\xop \xel_i^{-1}\xop \xel_{i\mp1} \\
  (\xel_j)\bel_i^{\pm1}&=&\xel_j{\rm\ if\ }i\ne j.
\end{eqnarray*}

We note the possible ways a generator or its inverse of the braid 
group can change the leading variable of an element of the free 
group written in reduced form:

\begin{lemma}\label{thatlemma}
If the reduced form of a word ${f}\in\xgp$ begins with $\xel_m$, 
then the reduced form of $({f})\bel_k^{\pm 1}$ begins 
with (this same) $\xel_m$ except in the following two cases (with the
displayed words in reduced form):
\begin{enumerate}
\item 
$(\xel_m\xop\xel_{m\pm1}^{-1}\xop\ldots)
  \bel_{m\pm1}^{\pm1} =
  \xel_{m\pm1}\xop\ldots$
\item
$(\xel_m\xop\ldots)\bel_m^{\pm1}=
   \xel_{m\pm1}\xop\xel_{m}^{-1}\xop\ldots$
\end{enumerate}
\end{lemma}

\bproof
Write $\bel$ for $\bel_k^{\pm1}$.

In each of the following cases we assume ${f}$ 
is reduced, that the braid action is applied to ${f}$, and then 
some unspecified but fixed reduction is applied to the word in 
$\xgp$ that results.

\begin{case}
$m\ne k$, and the leading $\xel_m$ is cancelled out in the reduction
of $({f})\bel$ by an $\xel_m^{-1}$ which is already present in the
reduced form of ${f}$.
\end{case}

This means that ${f}$ had the reduced form 
${f}=\xel_m\xop{f}_1\xop\xel_m^{-1}\xop{f}_2$.  Then $({f})\bel =
\xel_m\xop({f}_1)\bel\xop\xel_m^{-1}\xop({f}_2)\bel.$ But then
$({f}_1)\bel=\xid$, so ${f}_1=\xid$, and ${f}$ was not in reduced
form, contradiction.

\begin{case}
$m\ne k$, and the leading $\xel_m$ is cancelled out by an
$\xel_m^{-1}$ which was produced by the braid action on ${f}$.
\end{case}

An examination of cases shows that the 
only ways that an $\xel_m^{-1}$ can be so produced are:
\begin{enumerate}
\item
$(\xel_{m\pm1}^{-1})\bel_{m\pm1}^{\pm1}
 = \xel_m^{-1}\xop\xel_{m\pm1}\xop\xel_{m\pm2}^{-1}$

\item
$(\xel_{m\pm1}^{-1})\bel_{m\pm1}^{\mp1}
 = \xel_{m\pm2}^{-1}\xop\xel_{m\pm1}\xop\xel_m^{-1}$
\end{enumerate}

This means that ${f}$ had the reduced form 
${f}=\xel_m\xop{f}_1\xop\xel_{m\pm1}^{-1}\xop{f}_2$ and 
that $k=m\pm1$.  

\begin{subcase}
$\sigma=\bel_{m\pm1}^{\pm1}$.
\end{subcase}

The first of these two possibilities develops as:
$({f})\bel_{m\pm1}^{\pm1} = \xel_m\xop({f}_1)\bel_{m\pm1}^{\pm1}\xop
 (\xel_m^{-1}\xop\xel_{m\pm1}\xop\xel_{m\pm2}^{-1})
 \xop({f}_2)\bel_{m\pm1}^{\pm1}$.
In order that this $\xel_m^{-1}$ should be the one which cancels 
the leading $\xel_m$, it must be that $({f}_1)\bel=\xid$, and so 
${f}_1=\xid$.  Hence the reduced form of ${f}$ 
must be ${f}=\xel_m\xop\xel_{m\pm1}^{-1}\xop{f}_2$.

This is the lefthand side of the first possibility of the lemma.  
To finish this case we need to ensure that the 
$\xel_{m\pm1}$ produced above, which is the leading 
factor of $({f})\bel$ before reducing,
remains the leading factor after reducing.

If it were to cancel out in the reduction, then it must do so by 
cancelling with an $\xel_{m\pm1}^{-1}$ produced by the 
action of $\bel$ on an $\xel_{m\pm1}$ inside ${f}_2$.

Hence ${f}_2={f}_3\xop\xel_{m\pm1}\xop{f}_4$ in reduced form, so that
${f}=\xel_m\xop\xel_{m\pm1}^{-1} 
\xop{f}_3\xop\xel_{m\pm1}\xop{f}_4$.  Then $({f})\bel =
\xel_{m\pm1}\xop\xel_{m\pm2}^{-1} \xop({f}_3)\bel\xop
\xel_{m\pm2}\xop\xel_{m\pm1}^{-1}\xop\xel_m  \xop({f}_4)\bel.$ In
order that this $\xel_{m\pm1}^{-1}$ be the one to cancel the leading
$\xel_{m\pm1}$, it must be that $\xel_{m\pm2}^{-1}
\xop({f}_3)\bel\xop \xel_{m\pm2}=\xid$, and so ${f}_3=\xid$.  But
then ${f}$ was not in reduced form.

\begin{subcase}
$\sigma=\bel_{m\pm1}^{\mp1}$.
\end{subcase}

The second of these two possibilities develops as:
$({f})\bel_{m\pm1}^{\mp1} = \xel_m\xop({f}_1)\bel_{m\pm1}^{\mp1}\xop
(\xel_{m\pm2}^{-1}\xop\xel_{m\pm1}\xop\xel_m^{-1}) 
\xop({f}_2)\bel_{m\pm1}^{\mp1}.$ In order that this $\xel_m^{-1}$ be
the one to cancel the leading $\xel_m$, it must be that 
$\xel_{m\pm2}^{-1}
\xop({f}_3)\bel\xop
\xel_{m\pm2}=\xid$, and so
${f}_3=\xid$.
But then ${f}$ was not in reduced form.

\begin{case}
$k=m$.
\end{case}

Let ${f}=\xel_m\xop{f}_1$ in reduced form.
Then
$({f})\bel = 
\xel_{m\pm1}\xop\xel_m^{-1}\xop
\xel_{m\mp1}\xop({f}_1)\bel_m^{\pm1}.$

\begin{subcase}
The second factor of $(f)\bel$, $\xel_m^{-1}$, 
does not cancel in the reduction of $({f})\bel$.
\end{subcase}

If this $\xel_m^{-1}$ does not cancel in the reduction, 
then this case leads to the second possibility stated in the 
lemma.

\begin{subcase}
The second factor of $(f)\bel$, $\xel_m^{-1}$, 
does cancel in the reduction of 
$({f})\bel$.
\end{subcase}

In order that this $\xel_m^{-1}$ cancel, 
there must be an $\xel_m^{-1}$ in the 
reduced form of ${f}_1$, i.e., ${f}_1={f}_2\xop\xel_m^{-1}\xop{f}_3$ 
and so ${f}=\xel_m\xop{f}_2\xop\xel_m^{-1}\xop{f}_3$ in reduced form.

Hence
$({f})\bel = \xel_{m\pm1}\xop\xel_m^{-1}\xop\xel_{m\mp1}\xop
({f}_2)\bel\xop\xel_{m\mp1}^{-1}\xop\xel_m\xop
\xel_{m\pm1}^{-1}\xop({f}_3)\bel.$
In order that this $\xel_m$ be the one to cancel the second 
$\xel_m^{-1}$, it must be that
$\xel_{m\mp1}\xop({f}_2)\bel\xop\xel_{m\mp1}^{-1} 
=\xid$, and so ${f}_2=\xid$.
But then
${f}$ was not in reduced form.

This completes the proof of the lemma.
\eproof

\begin{definition}
Say that a reduced word in $\xgp$ leans right at $n$ if it begins 
with either $\xel_m$ or $\xel_n\xop\xel_m^{-1}$ with $m>n$.
\end{definition}

\begin{lemma}
If $f\in\xgp$ leans right at $n$, then so does $(f)\bel_i^{\pm1}$, 
where $\bel_i^{\pm1}\ne\bel_n^{-1}$.
\end{lemma}

\bproof
We first verify that if $f$ leans right at $n$, then the leading 
factor of the reduced form of $(f)\bel_i^{\pm1}$ is some $x_m$ with 
$m\ge n$.

By lemma \ref{thatlemma},
the index of the leading factor can 
change by at most $1$.

Hence assume by way of contradiction that $m=n$ and the action of 
$\bel_i^{\pm1}$ on $f$ leaves it with a leading factor of $\xel_{n-1}$.

The first possibility afforded by lemma \ref{thatlemma} is
\[(\xel_n\xop\xel_{n-1}^{-1}\xop\ldots)\bel_{n-1}^{-1}=
 \xel_{n-1}\xop\ldots\]
but this is ruled out by the requirement that the index of the 
second factor exceed $n$.

The second possibility is
$(\xel_n\xop\ldots)\bel_n^{-1}=
\xel_{n-1}\xop\xel_n^{-1}\xop\ldots$
and this is ruled out by $\bel_i^{\pm1}\ne\bel_n^{-1}$.

Hence if $f$ leans right at $n$ then $(f)\bel_i^{\pm1}$ leads with 
$\xel_n$ or $\xel_m$ with $m>n$.

Now we verify that if $f$ leans right at $n$, and $(f)\bel_i^{\pm1}$ 
starts with $x_n$, then its second factor is $x_m^{-1}$ with 
$m>n$.

Two possibilities arise.  Either $f$ starts with 
$\xel_{n+1}$ or with $\xel_n$.

In the first case, where $f$ starts with $\xel_{n+1}$, the lemma 
gives two possibilities where $(f)\bel_i^{\pm1}$ should start with 
$\xel_n$.

The first is $(\xel_{n+1}\xop\xel_n^{-1}\xop\ldots)\bel_n^{-1}=
\xel_n\xop\ldots$, but this is ruled out by
$\bel_i^{\pm1}\ne\bel_n^{-1}$.

The second is $(\xel_{n+1}\xop\ldots)\bel_{n+1}^{-1}=
\xel_n\xop\xel_{n+1}^{-1}\xop\ldots$, and this result (which is in 
reduced form) leans right at $n$.

In the second case, where $f$ starts with $\xel_n$, since $f$ 
leans right at $n$, $f$ must start with $\xel_n\xop\xel_m^{-1}$ with 
$m>n$.  We note that a dual version of the lemma shows that a 
reduced word starting with $\xel_m^{-1}$ and acted upon by a braid 
generator or its inverse either does not change its leading 
factor, or the leading factor becomes $\xel_{m\pm1}^{-1}$.  
Hence the only concern is $f=\xel_n\xop\xel_{n+1}^{-1}\xop\ldots$ 
acted on by a braid generator or its inverse that does not change 
the leading $\xel_n$ and does change $\xel_{n+1}^{-1}\xop\ldots$ 
into $\xel_n^{-1}\xop\ldots$, which contradicts that the leading 
$\xel_n$ is not cancelled.

Therefore if $f$ leans right at $n$, and $(f)\bel_i^{\pm1}$ leads with 
$\xel_n$, then its second factor is $\xel_m^{-1}$ with $m>n$.

This completes the proof.
\eproof

\bproof
To prove the $\sigma_n$-proposition, 
assume that ${p}\in\bgp$ is formed from a product of generators of 
$\bgp$ and their inverses, with at least one $\bel_n$ and no 
$\bel_n^{-1}$.

Write ${p}={p}_1\bop\bel_n\bop{p}_2$ 
where $\bel_n^{\pm1}$ doesn't occur in $p_1$ and 
$\bel_n^{-1}$ doesn't occur in ${p}_2$.

Then
\begin{eqnarray*}
 (\xel_n){p} &=& (\xel_n){p}_1\bop\bel_n\bop{p}_2 \\
  &=& (\xel_n)\bel_n\bop{p}_2 \\
  &=& (\xel_{n+1}\xop\xel_n^{-1}\xop\xel_{n-1}){p}_2
\end{eqnarray*}

But $\xel_{n+1}\xop\xel_n^{-1}\xop\xel_{n-1}$ leans right at $n$, 
and hence so does
$(\xel_{n+1}\xop\xel_n^{-1}\xop\xel_{n-1}){p}_2$.  Since $\xel_n$ 
does not lean right at $n$, they cannot be equal, and so 
${p}\ne\bid$.
\eproof


\begin{thebibliography}{Xxx 99x}

\bibitem[Bir 75]{bir:bg}{\sc J.~Birman}, Braids, Links, and Mapping
Class Groups, {\em Annals of Mathematics Studies\/} {\bf 82} (1975).

\bibitem[Deh 89a]{deh:fdg}{\sc P.~Dehornoy}, Free Distributive
Groupoids, {\em Journal of Pure and Applied Algebra\/} {\bf 61} 
(1989), 123--146.

\bibitem[Deh 89b]{deh:gl}{\sc P.~Dehornoy}, Sur la structure des
ger\-bes libres, {\em Comp\-tes-Rendus de l'Acad.~des Sciences de Paris\/}
{\bf 309-I} (1989), 143--148.

\bibitem[Deh 92]{deh:bg}{\sc P.~Dehornoy}, Braid Groups and Left
Distributive Structures, {\em preprint\/} (1992).

\bibitem[Lav 92]{lav:ld}{\sc R.~Laver}, The Left Distributive Law and
the Freeness of an Algebra of Elementary Embeddings, {\em Advances in
Mathematics\/} {\bf 91} (1992), 209--231.
\end{thebibliography}
\end{document}